\documentclass[12pt]{amsart}
\usepackage{latexsym}
\usepackage{amscd,amssymb,times,fullpage}
\newtheorem{thm}{Theorem}
\newtheorem{prop}{Proposition}
\newtheorem{lem}{Lemma}
\newtheorem{cor}{Corollary}

\theoremstyle{remark}
\newtheorem{rem}{Remark}
\newtheorem{ex}{Example}

\theoremstyle{definition}

\newcommand{\bZ}{\mathbb Z}
\newcommand{\bR}{\mathbb R}
\newcommand{\bC}{\mathbb C}
\newcommand{\bP}{\mathbb P}

\newcommand{\blowup}{\overline{\bC P^2}}

\newcommand{\JJ}{\mathcal J}

\newcommand{\C}{\mathbb{ C}}

\newcommand{\lra}{\longrightarrow}

\newcommand{\OO}{\mathcal{ O}}

\newcommand{\ra}{{\rightarrow}}

\newcommand{\Z}{\mathbb{ Z}}

\newcommand{\Spc}{\mathrm{Spin}^c}

\newcommand{\fs}{\mathfrak{ s}}

\title{Minimality and irreducibility of symplectic four-manifolds}
\author{M.~J.~D.~Hamilton}
\address{Mathematisches Institut, Ludwig-Maximilians-Universit\"at M\"unchen,
Theresienstr.~39, 80333 M\"unchen, Germany}
\email{Mark.Hamilton@mathematik.uni-muenchen.de}
\author{D.~Kotschick}
\address{Mathematisches Institut, Ludwig-Maximilians-Universit\"at M\"unchen,
Theresienstr.~39, 80333 M\"unchen, Germany}
\email{dieter@member.ams.org}
\thanks{The second author wishes to thank Y.~Eliashberg for his hospitality at Stanford University
while some of this work was carried out.}
\date{September 5, 2005; MSC 2000: primary 57R17, 57R57, 32J15, secondary 53D35, 32G05, 32J27}

\begin{document}

\begin{abstract}
    We prove that all minimal symplectic four-manifolds are essentially irreducible. 
    We also clarify the relationship between holomorphic and symplectic minimality 
    of K\"ahler surfaces. This leads to a new proof of the deformation-invariance of
    holomorphic minimality for complex surfaces with even first Betti number which 
    are not Hirzebruch surfaces.
\end{abstract}

\maketitle

\section{Introduction and statement of results}

In this paper we discuss certain geometric and topological properties of 
symplectic four-manifolds. Our main concern is the notion of minimality, 
and its topological consequences. We shall extend to manifolds with 
$b_{2}^{+}=1$ the irreducibility result proved in~\cite{irred,Bourbaki} for  
the case that $b_{2}^{+}>1$. We also show that holomorphic and symplectic 
minimality are equivalent precisely for those K\"ahler surfaces which are not 
Hirzebruch surfaces.
Together with work of Buchdahl~\cite{B}, this yields a new proof of the 
deformation-invariance of holomorphic minimality for complex surfaces
with even first Betti number, again with the exception of Hirzebruch surfaces.

\subsection{Minimality}\label{ss:mini}

A complex surface is said to be minimal if it contains no holomorphic sphere 
of selfintersection $-1$, see for example~\cite{BPV}. A symplectic four-manifold 
is usually considered to be minimal if it contains no symplectically embedded 
sphere of selfintersection $-1$, see for example~\cite{McD1,Go}. In the case 
of a K\"ahler surface both notions of minimality can be considered, but it is 
not at all obvious whether they agree. In the recent literature on symplectic 
four-manifolds there are frequent references to (symplectic) minimality, and 
often K\"ahler surfaces are considered as examples, but we have found no explicit 
discussion of the relationship between the two definitions in print, compare 
e.~g.~\cite{McD1,McD2,McD3,R,Go,Bourbaki,GoS}.

An embedded holomorphic curve in a K\"ahler manifold is a symplectic 
submanifold. Therefore, for K\"ahler surfaces symplectic minimality implies 
holomorphic minimality. The following counterexample to the converse 
should be well known:
\begin{ex}\label{ex:H}
    Let $X_{n}=\bP (\OO\oplus\OO(n))$ be the $n^{\rm{th}}$ Hirzebruch surface. 
    If $n$ is odd and $n>1$, then $X_{n}$ is holomorphically minimal 
    but not symplectically minimal.
    \end{ex}
In Section~\ref{s:mini} below we explain this example in detail, and then we 
prove that there are no other counterexamples:
\begin{thm}\label{Hirzebruch} 
    A K\"ahler surface that is not a Hirzebruch surface $X_n$ with 
    $n$ odd and $n>1$ is holomorphically minimal if and only if it is 
    symplectically minimal. 
\end{thm}
A proof can be given using the known calculations of Seiberg--Witten 
invariants of K\"ahler surfaces. Using Seiberg--Witten theory, it turns 
out that for non-ruled K\"ahler surfaces symplectic and holomorphic 
minimality coincide because they are both equivalent to smooth minimality, 
that is, the absence of smoothly embedded $(-1)$-spheres. The case of 
irrational ruled surfaces is elementary.

Such a proof is not satisfying conceptually, because the basic notions of 
symplectic topology should be well-defined without appeal to results in 
gauge theory. Therefore, in Section~\ref{s:mini} we give a proof of 
Theorem~\ref{Hirzebruch} within the framework of symplectic topology, using 
Gromov's theory of $J$-holomorphic curves. We shall use results of 
McDuff~\cite{McD1} for which Gromov's compactness theorem is crucial. 
Essentially the same argument can be used to show that symplectic minimality 
is a deformation-invariant property, see Theorem~\ref{t:def}. This natural result 
is lurking under the surface of McDuff's papers~\cite{McD1,McD2,McD3}, and 
is made explicit in~\cite{MP}, compare also~\cite{Rjdg,R}. Of course this result 
is also a corollary of Taubes's deep work in~\cite{swg,main,Taubes,Bourbaki}, 
where he showed, among other things, that if there is a smoothly embedded 
$(-1)$-sphere, then there is also a symplectically embedded one. 

In Section~\ref{s:mini} we shall also prove that for compact complex surfaces with even 
first Betti number which are not Hirzebruch surfaces holomorphic minimality is 
preserved under deformations of the complex structure. This result is known, and
is traditionally proved using the Kodaira classification, cf.~\cite{BPV}. The proof we 
give is intrinsic and independent of the classification. Instead, we combine the result 
of Buchdahl~\cite{B} with the deformation invariance of symplectic minimality and 
Theorem~\ref{Hirzebruch}.

\subsection{Irreducibility}
Recall that an embedded $(-1)$-sphere in a four-manifold gives rise to a connected sum 
decomposition where one of the summands is a copy of $\blowup$. For symplectic
manifolds no other non-trivial decompositions are known. Gompf~\cite{Go} conjectured 
that minimal symplectic four-manifolds are irreducible, meaning that in any smooth 
connected sum decomposition one of the summands has to be a homotopy sphere. 
In Section~\ref{s:irred} below we shall prove the following result in this direction:
\begin{thm}\label{sum}
    Let $X$ be a minimal symplectic $4$-manifold with $b_2^+=1$. If $X$ splits 
    as a smooth connected sum $X = X_1 \# X_2$, then one of the $X_i$ is an 
    integral homology sphere whose fundamental group has no non-trivial finite 
    quotient. 
\end{thm}
For manifolds with $b_2^+ >1$ the corresponding result was first proved 
in~\cite{irred} and published in~\cite{Bourbaki}. As an immediate consequence 
of these results we verify Gompf's irreducibility conjecture in many cases:
\begin{cor} 
    Minimal symplectic $4$-manifolds with residually finite fundamental groups 
    are irreducible.
\end{cor}

To prove Theorem~\ref{sum} we shall follow the strategy of the proof for 
$b_{2}^{+}>1$ in~\cite{irred,Bourbaki}. In particular we shall use the deep 
work of Taubes~\cite{swg,main,Taubes}, which produces symplectic submanifolds 
from information about Seiberg--Witten invariants. What is different in the case 
$b_{2}^{+}=1$, is that the Seiberg--Witten invariants depend on chambers, and 
one has to keep track of the chambers one is working in. 

In addition to conjecturing the irreducibility of minimal symplectic four-manifolds, 
Gompf~\cite{Go} also raised the question whether minimal non-ruled symplectic 
four-manifolds satisfy $K^{2}\geq 0$, where $K$ is the canonical class. For 
manifolds with $b_{2}^{+}>1$ this was proved by Taubes~\cite{swg,main}, compare 
also~\cite{Bourbaki,Taubes}. The case $b_{2}^{+}=1$ was then treated by Liu~\cite{Liu}, 
who refers to this question as ``Gompf's conjecture''. Liu~\cite{Liu} also proved that 
minimal symplectic four-manifolds which are not rational or ruled satisfy 
$K\cdot\omega\geq 0$. We shall use Liu's inequalities to keep track of the chambers 
in our argument. Although the results of Liu~\cite{Liu}, and also those of Li--Liu~\cite{LL2,LL3}, 
are related to Theorem~\ref{sum}, this theorem does not appear there, or anywhere 
else in the literature that we are aware of.

\section{Notions of minimality}\label{s:mini}

First we discuss the Hirzebruch surfaces $X_{n}=\bP (\OO\oplus\OO(n))$, with $n$ odd and
$>1$, in order to justify the assertions made in Example~\ref{ex:H} in the Introduction.

If $n=2k+1$, consider the union of a holomorphic section $S$ of $X_n$ of selfintersection 
$-n$ and of $k$ disjoint parallel copies of the fiber $F$. This reducible holomorphic curve 
can be turned into a symplectically embedded sphere $E$ by replacing each of the transverse 
intersections of $S$ and $F$ by a symplectically embedded annulus. Then
$$
E\cdot E = (S+kF)^{2} = S\cdot S + 2k \  S\cdot F = -n +2k = -1 \ .
$$
This shows that $X_{n}$ is not symplectically minimal. To see that it is holomorphically
minimal, note that a homology class $E$ containing a smooth holomorphic
$(-1)$-sphere would satisfy $E^2=K\cdot E=-1$, and would therefore be $S+kF$, as 
above. However, this class has intersection number 
$$
E\cdot S = (S+kF)\cdot S = -n+k=-k-1<0
$$
with the smooth irreducible holomorphic curve $S$. Therefore, $E$ can only contain 
a smooth irreducible holomorphic curve if $E=S$, in which case $k=0$ and $n=1$.

Next we prove that for all other K\"ahler surfaces symplectic and holomorphic minimality
are equivalent.
\begin{proof}[Proof of Theorem~\ref{Hirzebruch}]
In view of the discussion in~\ref{ss:mini} above, we only have to prove that if $(X,\omega)$ 
is a K\"ahler surface which is not a Hirzebruch surface $X_n$ with $n$ odd and $n>1$, 
then holomorphic minimality implies symplectic minimality.

We start by assuming that $(X,\omega)$ is not symplectically minimal, 
so that it contains a smoothly embedded $(-1)$-sphere $E\subset X$ with 
$\omega\vert_{E}\neq 0$. Orient $E$ so that $\omega\vert_{E}>0$, and 
denote by $[E]\in H_{2}(X;\Z)$ the corresponding homology class. The 
almost complex structures $J$ compatible with $\omega$ are all homotopic 
to the given integrable $J_{\infty}$; in particular their canonical 
classes agree with the canonical class $K$ of the K\"ahler structure. It 
is elementary to find a compatible $J$ for which the sphere $E$ with the 
chosen orientation is $J$-holomorphic. Therefore $E$ satisfies the 
adjunction formula
$$
g(E)=1+\frac{1}{2}(E^{2}+K\cdot E) \ .
$$
We conclude that $K\cdot E=-1$. (Note that the orientation of $E$ is 
essential here.) This implies in particular that the expected 
dimension of the moduli space of $J$-holomorphic curves in the 
homology class $[E]$ vanishes.

Let $\JJ$ be the completion--with respect to a suitable Sobolov 
norm--of the space of $C^{\infty}$ almost complex structures compatible 
with $\omega$, cf.~\cite{MS}. McDuff has proved that, for almost complex 
structures $J$ from an everywhere dense subset in $\JJ$, there is a 
unique smooth $J$-holomorphic sphere $C$ in the homology class $[E]$, see 
Lemma~3.1 in~\cite{McD1}. 

The uniqueness implies that the curve $C$ varies smoothly with $J$. One 
then uses Gromov's compactness theorem for a smooth family of almost complex 
structures to conclude that for all $J$, not necessarily generic, there 
is a unique $J$-holomorphic representative of the homology class $[E]$ 
which, if it is not a smooth curve, is a reducible curve $C=\sum_{i}C_{i}$ 
such that each $C_{i}$ is a smooth $J$-holomorphic sphere. Compare again 
Lemma~3.1 in~\cite{McD1} and~\cite{MS}. (In these references reducible 
$J$-holomorphic curves are called cusp curves.)

Let $J_{j}$ be a sequence of generic almost complex structures in $\JJ$ 
which converges to the integrable $J_{\infty}$ as $j\ra\infty$. For each 
$J_{j}$ there is a smooth $J_{j}$-holomorphic sphere $E_{j}$ in the homology 
class $[E]$. As $j\ra\infty$, the $E_{j}$ converge 
weakly to a possibly reducible $J_{\infty}$-holomorphic curve 
$E_{\infty}$. If $E_{\infty}$ is irreducible, then it is a holomorphic 
$(-1)$-sphere, showing that $(X,J_{\infty})$ is not holomorphically minimal. 
If $E_{\infty}$ is reducible, let 
$$
E_{\infty}=\sum_{i=1}^{k}m_iC_{i}
$$
be the decomposition into irreducible components. The multiplicities $m_i$ are positive 
integers. Each $C_{i}$ is an embedded sphere, and therefore the adjunction formula implies
$$
C_{i}^{2}+K\cdot C_{i}=-2 \ .
$$
Multiplying by $m_i$ and summing over $i$ we obtain
$$
\sum_{i=1}^{k}m_iC_{i}^{2}+K\cdot\sum_{i=1}^{k}m_iC_{i}=-2\sum_{i=1}^{k}m_i \ .
$$
Now the second term on the left hand side equals $K\cdot E=-1$, so 
that we have
$$
\sum_{i=1}^{k}m_iC_{i}^{2}=1-2\sum_{i=1}^{k}m_i \ .
$$
It follows that there is an index $i$ such that $C_{i}^{2}\geq -1$. 
If $C_{i}^{2}=-1$ for some $i$, then we again conclude that $(X,J_{\infty})$ 
is not holomorphically minimal. If $C_{i}^{2}\geq 0$ for some $i$, 
then $(X,J_{\infty})$ is birationally ruled or is rational, 
cf.~Proposition~4.3 in Chapter~V of~\cite{BPV}. Thus, if 
it is holomorphically minimal, it is either a minimal ruled surface 
or $\C P^{2}$, but the latter is excluded by our assumption that 
$(X,J_{\infty})$ is not symplectically minimal. If $(X,J_{\infty})$ 
were ruled over a surface of positive genus, $X\stackrel{\pi}{\lra}B$, 
then the embedding of the $(-1)$-sphere $E$ would be homotopic to a 
map with image in a fiber, because $\pi\vert_{E}\colon E\ra B$ would 
be homotopic to a constant. But this would contradict the fact that 
$E$ has non-zero selfintersection.

Thus we finally reach the conclusion that $(X,J_{\infty})$ is ruled 
over $\C P^{1}$. If it is holomorphically minimal, then it is a 
Hirzebruch surface $X_n$ with $n$ odd and $n>1$, because $X_{1}$ is 
not holomorphically minimal, and $X_{2k}$ has even intersection form 
and is therefore symplectically minimal.

This concludes the proof of Theorem~\ref{Hirzebruch}.
\end{proof}
\begin{rem}
    We have used the fact that the existence of a rational holomorphic curve of 
    non-negative selfintersection in a complex surface implies that 
    the surface is rational or ruled. Such a statement also holds in 
    the symplectic category, cf.~\cite{McD1}, but we do not need that 
    here.
    \end{rem}
    
    The exposition of the proof of Theorem~\ref{Hirzebruch} can be 
    shortened considerably if one simply uses McDuff's Lemma~3.1 
    from~\cite{McD1} as a black box. We have chosen to include some of 
    the details so that the reader can see that the degeneration of 
    the $J_{j}$-holomorphic curves $E_{j}$ as $j\ra\infty$ is the 
    exact inverse of the regeneration used in the discussion of 
    Example~\ref{ex:H}. 
    
    The following theorem, Proposition~2.3.A in~\cite{MP}, can be proved by 
    essentially the same argument, allowing the symplectic form to vary smoothly,
    compare also~\cite{McD1,R}:
\begin{thm}[\cite{MP}]\label{t:def}
    Symplectic minimality is a deformation-invariant property of 
    compact symplectic four-manifolds.
    \end{thm}

   Note that holomorphic minimality of complex surfaces is not invariant 
    under deformations of the complex structure. In the K\"ahler case the 
    Hirzebruch surfaces $X_{n}$ with $n$ odd are all deformation-equivalent, 
    but are non-minimal for $n=1$ and minimal for $n>1$. In the non-K\"ahler 
    case there are other examples among the so-called class $VII$ surfaces.

    For complex surfaces of non-negative Kodaira dimension it is true that holomorphic 
    minimality is deformation-invariant, but the traditional proofs for this are 
    exceedingly cumbersome, see for example~\cite{BPV}, section~7 of Chapter~VI, 
    where it is deduced from the Kodaira classification and a whole array of additional 
    results. For the case of even first Betti number we now give a direct proof, which 
    does not use the classification.
 \begin{thm}
Let $X$ be a holomorphically minimal compact complex surface with even first 
Betti number, which is not a Hirzebruch surface $X_n$ with $n$ odd. Then any 
surface deformation equivalent to $X$ is also holomorphically minimal.
 \end{thm}
 \begin{proof}
Let $X_t$ with $t\in [0,1]$ be a smoothly varying family of complex surfaces such 
that $X_0=X$.  Buchdahl~\cite{B} has proved that every compact complex surface 
with even first Betti number is K\"ahlerian, without appealing to any classification 
results. Thus, each $X_t$ is K\"ahlerian, and we would like to choose K\"ahler 
forms $\omega_0$ and $\omega_1$ on $X_0$ and $X_1$ respectively, which can 
be joined by a smooth family of symplectic forms $\omega_t$. There are two ways 
to see that this is possible. 

On the one hand, Buchdahl~\cite{B} characterizes the K\"ahler classes, and one can 
check that one can choose a smoothly varying family of K\"ahler classes for $X_t$, 
which can then be realized by a smoothly varying family of K\"ahler metrics. On the 
other hand, we could just apply Buchdahl's result for each value of the parameter $t$ 
separately, without worrying about smooth variation of the K\"ahler form with the parameter,
and then construct a smooth family $\omega_t$ of symplectic not necessarily K\"ahler 
forms from this, cf.~\cite{R} Proposition 2.1. In detail, start with arbitrary K\"ahler forms
$\omega_t$ on $X_t$. As the complex structure depends smoothly on $t$, there 
is an open neighbourhood of each $t_0\in [0,1]$ such that $\omega_{t_0}$ is a compatible 
symplectic form for all $X_s$ with $s$ in this neighbourhood of $t_0$. By compactness 
of $[0,1]$, we only need finitely many such open sets to cover $[0,1]$. On the 
overlaps we can deform these forms by linear interpolation, because the space of 
compatible symplectic forms is convex. In this way we obtain a smoothly varying family
of symplectic forms.

Now $X=X_0$ was assumed to be holomorphically minimal and not a Hirzebruch surface
$X_n$ with odd $n$. Therefore, Theorem~\ref{Hirzebruch} shows that $X_0$ is 
symplectically minimal, and Theorem~\ref{t:def} then implies that $X_1$ is also 
symplectically minimal. The easy direction of Theorem~\ref{Hirzebruch} shows that
$X_1$ is holomorphically minimal.
 \end{proof}
 Let us stress once more that this result is not new, but its proof is. The above proof does 
 not use the Kodaira classification. The only result we have used from the traditional theory  
 of complex surfaces is that a surface containing a holomorphic sphere of positive  square 
 is rational, which entered in the proof of Theorem~\ref{Hirzebruch}. We have not used the 
 generalization of this result to symplectic manifolds, and we have not used any Seiberg--Witten 
 theory either. Our proof does depend in an essential way on the work of Buchdahl~\cite{B}. 
 Until that work, the proof that complex surfaces with even first Betti numbers are K\"ahlerian 
 depended on the Kodaira classification.
      
\section{Connected sum decompositions of minimal symplectic 
four-manifolds}\label{s:irred}

In this section we prove restrictions on the possible connected sum decompositions of a 
minimal symplectic four-manifold with $b_2^+=1$, leading to a proof of Theorem~\ref{sum}.
To do this we have to leave the realm of symplectic topology and use Seiberg--Witten gauge theory.

Let $X$ be a closed oriented smooth $4$-manifold with $b_2^+(X)=1$. We fix a $\Spc$ structure 
$\fs$ and a metric $g$ on $X$ and consider the Seiberg--Witten equations for a positive spinor $\phi$
and a $\Spc$ connection $A$:
\begin{eqnarray*}
D^+_A \phi & = & 0 \\
F^+_{\hat A} & = & \sigma(\phi,\phi)+\eta \ ,
\end{eqnarray*} 
where the parameter $\eta$ is an imaginary-valued $g$-self-dual $2$-form. Here $\hat A$ denotes
the $U(1)$-connection on the determinant line bundle induced from $A$, so that $F_{\hat A}$
is an imaginary-valued $2$-form.
A reducible solution of the Seiberg--Witten equations 
is a solution with $\phi=0$. 

For every Riemannian metric $g$ there exists a $g$-self-dual harmonic $2$-form $\omega_g$ 
with $[\omega_g]^2=1$. Because $b_2^+(X)=1$, this $2$-form is determined by $g$ up to a sign. 
We choose a forward cone, i.~e.~one of the two connected components of 
$\{\alpha \in H^2(X;\mathbb{R}) \ \vert \ \alpha^2 > 0\}$. 
Then we fix $\omega_g$ by taking the form whose cohomology class lies 
in the forward cone.

Let $L$ be the determinant line bundle of the $\Spc$ structure $\fs$. The curvature 
$F_{\hat A}$ represents $\frac{2\pi}{i} c_1(L)$ in cohomology, and every form which represents this class can be 
realized as the curvature of  $\hat A$ for a $\Spc$ connection $A$. 
For given $(g,\eta)$ there exists a reducible solution of the Seiberg--Witten equations if and only if there is a
$\Spc$ connection $A$ such that  $F_{\hat A}^+=\eta$, equivalently
$(c_1(L)-\frac{i}{2\pi}\eta) \cdot \omega_g=0$. 
Define the discriminant of the parameters $(g,\eta)$ by 
\begin{equation*}
\Delta_L(g,\eta)=(c_1(L)-\frac{i}{2\pi}\eta) \cdot \omega_g \ .
\end{equation*}
One divides the space of parameters $(g,\eta)$ for which there are no reducible solutions into the plus and minus 
chambers according to the sign of the discriminant. Two pairs of parameters $(g_1,\eta_1)$ and $(g_2,\eta_2)$ 
can be connected by a path avoiding reducible solutions if and only if their discriminants have the same sign, 
i.~e.~if and only if they lie in the same chamber. A cobordism argument then shows that the Seiberg--Witten 
invariant is the same for all parameters in the same chamber. In this way we get the invariants $SW_+(X,\fs)$,
$SW_-(X,\fs)$ which are constant on the corresponding chambers. 

Suppose now that $X$ has a symplectic structure $\omega$. Then $\omega$ determines an orientation of $X$ and 
a forward cone in $H^2(X;\mathbb{R})$. We will take the chambers with respect to this choice. Moreover, $\omega$ 
determines a canonical class $K$ and a $\Spc$ structure $\fs_{K^{-1}}$ with determinant $K^{-1}$. One can obtain 
every other $\Spc$ structure by twisting $\fs_{K^{-1}}$ with a line bundle $E$, to obtain $\fs_{K^{-1}} \otimes E$. 
This $\Spc$ structure has determinant $K^{-1}\otimes E^2$. 

The {\it Taubes chamber} is the chamber determined by parameters $(g,\eta)$ with  
$g$ chosen such that it is almost K\"ahler with $\omega_g=\omega$ and 
\begin{equation*}
\eta = F_{\hat A_0}^+-\frac{i}{4}r\omega \quad\mbox{with}\quad r \gg 0 \ ,
\end{equation*}
where $\hat A_0$ is a canonical connection on $K^{-1}$. We have the following:
\begin{lem} 
The Taubes chamber is the minus chamber, for the choice of forward cone as above.
\end{lem}
\begin{proof} We have 
\begin{eqnarray*} (c_1(-K)-\frac{i}{2\pi}\eta) \cdot \omega_g & = & (\frac{i}{2\pi}F_{\hat{A_0}}-\frac{i}{2\pi}F_{\hat{A_0}}^+-\frac{1}{8\pi}r\omega)\cdot \omega \\
           & = & (\frac{i}{2\pi}F_{\hat{A_0}}^--\frac{1}{8\pi}r\omega)\cdot \omega \\
           & = & -\frac{1}{8\pi}r\omega^2 < 0 \ ,
\end{eqnarray*}
because the wedge product of a self-dual and an anti-self-dual two-form vanishes.
\end{proof}           
The following theorem is due to Taubes~\cite{T1,swg,main}, compare~\cite{LL2,LL3} for the case $b_2^+=1$.  
\begin{thm}\label{Taubes}
The Seiberg--Witten invariant in the minus chamber for the canonical $\Spc$ structure is non-zero. 
More precisely, $SW_-(X,\fs_{K^{-1}})=\pm 1$. Moreover, if $SW_-(X,\fs_{K^{-1}} \otimes E)$ is non-zero and $E \neq 0$, 
then for a generic $\omega$-compatible almost complex structure $J$, the Poincar\'{e} dual of the Chern class of $E$ 
can be represented by a smooth $J$-holomorphic curve $\Sigma \subset X$. 
\end{thm}

We have the following more precise version of the second part of Theorem~\ref{Taubes}, which is also due to Taubes.  
\begin{prop}\label{marcolli} 
Suppose $SW_-(X,\fs_{K^{-1}}\otimes E)$ is non-zero, and $E \neq 0$. Then for a generic almost complex structure 
$J$ compatible with $\omega$ there exist disjoint embedded $J$-holomorphic curves $C_i$ in $X$ such that 
 \begin{equation*}
PD(c_1(E))=\sum_{i=1}^nm_i[C_i],
\end{equation*}
where each $C_i$ satisfies $K\cdot C_i\leq C_i\cdot C_i$ and each multiplicity $m_i$ is equal to $1$, except possibly 
for those $i$ for which $C_i$ is a torus with self-intersection zero.
\end{prop}
This depends on a transversality result for $J$-holomorphic curves, see Proposition~7.1 in~\cite{main} and also~\cite{Taubes,Bourbaki}.
Proposition~\ref{marcolli} immediately implies the following:
\begin{cor}\label{minimal} 
If $SW_-(X,\fs_{K^{-1}}\otimes E) \neq 0$ with $E^2 <0$ then $X$ contains an embedded 
symplectic $(-1)$-sphere $\Sigma$.
\end{cor}
\begin{proof} 
Choose a generic compatible almost complex structure $J$ as in Proposition~\ref{marcolli}, and consider $E=\sum_i m_iC_i$. 
Then $E^2 = \sum_i m_i^2C_i^2$ because the $C_i$ are disjoint, hence $C_j^2 < 0$ for some $j$. We can compute 
the genus of $C_j$ from the adjunction formula:
\begin{equation*}
g(C_j)=1+\frac{1}{2}(C_j\cdot C_j+K\cdot C_j)\leq 1+ C_j\cdot C_j \leq 0.
\end{equation*}
Hence $\Sigma=C_j$ is a sphere with self-intersection number $-1$. 
\end{proof}
After these preparations we can now prove Theorem~\ref{sum}.

\begin{proof}[Proof of Theorem~\ref{sum}]
Let $(X,\omega)$ be a closed symplectic $4$-manifold with $b_2^+=1$. We denote by $K$ both the first Chern class
of any compatible almost complex structure, and the complex line bundle with this Chern class.

First, suppose that $(X,\omega)$ is symplectically minimal and rational or ruled. Then, by the classification of ruled symplectic 
four-manifolds, $X$ is diffeomorphic either to $\bC P^2$, to an even Hirzebruch surface, or to a geometrically ruled K\"ahler 
surface over a complex curve of positive genus, compare e.~g.~\cite{MS}. These manifolds are all irreducible for purely 
topological reasons. This is clear for $\bC P^2$ and for the even Hirzebruch surfaces, because the latter are diffeomorphic to 
$S^2\times S^2$. For the irrational ruled surfaces note that the fundamental group is indecomposable as a free product. 
Therefore, in any connected sum decomposition one of the summands is simply connected. If this summand were not a 
homotopy sphere, then the other summand would be a smooth four-manifold with the same fundamental group but with 
strictly smaller Euler characteristic than the ruled surface. This is impossible, because the irrational ruled surfaces realize 
the smallest possible Euler characteristic for their fundamental groups, compare~\cite{tani}.

Thus, we may assume that $(X,\omega)$ is not only symplectically minimal, but also not rational or ruled. 
Then Liu's results in~\cite{Liu} tell us that $K^2\geq 0$ and $K\cdot\omega\geq 0$.

If $X$ decomposes as a connected sum $X = M \# N$ then one of the summands, say $N$, has negative definite 
intersection form. Moreover, the fundamental group of $N$ has no non-trivial finite quotients, by Proposition~1 
of~\cite{KMT}. In particular $H_1(N;\mathbb{Z})=0$, and hence the homology and cohomology of $N$ are torsion-free. 
If $N$ is an integral homology sphere, then there is nothing more to prove.

Suppose $N$ is not an integral homology sphere. By Donaldson's theorem~\cite{D}, the intersection 
form of $N$ is diagonalizable over $\mathbb{Z}$. Thus there is a basis $e_1,...,e_n$ of $H^2(N;\mathbb{Z})$ 
consisting of elements with square $-1$ which are pairwise orthogonal. Write
$$
K = K_M+\sum_{i=1}^{n}a_ie_i \ ,
$$
with $K_M\in H^2(M;\bZ)$. The $a_i\in\bZ$ are odd, because $K$ is a characteristic vector. This shows in particular 
that $K$ is not a torsion class. Its orthogonal complement $K^\perp$ in $H^2(X;\bR)$ is then a hyperplane.
As $K^2\geq 0$ and $b_2^+(X)=1$, the hyperplane $K^\perp$ does not meet the positive cone. Thus 
Liu's inequality $K\cdot\omega\geq 0$ must be strict: $K\cdot\omega>0$.

Now we know $SW_-(X,\fs_{K^{-1}})=\pm 1$ from Taubes's result, where $\fs_{K^{-1}}$ is the $\Spc$ structure with 
determinant $K^{-1}$ induced by the symplectic form $\omega$. The inequality $(-K)\cdot\omega<0$ shows that 
a pair $(g,0)$ is in the negative, i.~e.~the Taubes chamber, whenever $g$ is almost K\"ahler with fundamental 
two-form $\omega$. As $K^\perp$ does not meet the positive cone, all pairs $(g,0)$ are in the negative chamber,
for all Riemannian metrics $g$.
We choose a family of Riemannian metrics $g_r$ on $X$ which pinches the neck connecting 
$M$ and $N$ down to a point as $r \rightarrow \infty$. For $r$ large we may assume that $g_r$ converges to 
metrics on the (punctured) $M$ and $N$, which we denote by $g_M$ and $g_N$. 
\begin{lem} \label{l}
If we choose the forward cone for $M$ to be such that it induces on $X$ the forward cone determined by the 
symplectic structure, then for every Riemannian metric $g'$ on $M$, the point $(g',0)$ is in the negative chamber 
of $M$ with respect to the $\Spc$ structure $\fs_M$ on $M$ obtained by restriction of $\fs_{K^{-1}}$.
\end{lem}
\begin{proof} 
The chamber is determined by the sign of $c_1(\fs) \cdot \omega_g$. 
We have 
\begin{eqnarray*}
0>(-K)\cdot\omega_{g_r}=c_1(\fs_{K^{-1}}) \cdot \omega_{g_r} & = & c_1(\fs_M) \cdot \omega_{g_r} + c_1(\fs_N) \cdot \omega_{g_r} \\
                          & \longrightarrow & c_1(\fs_M) \cdot \omega_{g_M} + c_1(\fs_N) \cdot \omega_{g_N},\,\mbox{as}\,r\rightarrow\infty \ .
\end{eqnarray*}
We know that $\omega_{g_N}$ is self-dual harmonic with respect to $g_N$,
and hence vanishes because $b_2^+(N)=0$. This implies that $c_1(\fs_{K^{-1}}) \cdot \omega_{g_r}$
converges to $c_1(\fs_M) \cdot \omega_{g_M}$ for $r\rightarrow\infty$. Thus
$$
c_1(\fs_M) \cdot \omega_{g_M}\leq 0 \ .
$$
However, we have $c_1(\fs_M) = K_M^{-1}$, and 
$$
K_M^{2}=K^2+\sum_{i=1}^{n}a_i^2\geq K^2+n\geq n\geq 1 \ ,
$$
showing that $K_M^\perp$ does not meet the positive cone of $M$. Thus $c_1(\fs_M) \cdot \omega_{g_M}<0$. Again
because $K_M^\perp$ does not meet the positive cone of $M$, this inequality holds for all metrics $g'$ on $M$.
\end{proof}
The degeneration of the $g_r$ as $r$ goes to infinity takes place in the negative chamber for $\fs_{K^{-1}}$, where the 
Seiberg--Witten invariant is $\pm 1$, and, by Lemma~\ref{l}, $g_M$ is in the negative chamber for $\fs_M$. 
It follows that $SW_-(M,\fs_M)=\pm 1$.

We now reverse the metric degeneration, but use a different $\Spc$ structure on $N$. Instead of using $\fs_N$ with
$c_1(\fs_N)=-\sum_{i=1}^na_ie_i$, we use the unique $\Spc$ structure $\fs_N'$ with $c_1(\fs_N')=a_1e_1-\sum_{i=2}^na_ie_i$.
For every metric on $N$ there is a unique reducible solution of the Seiberg--Witten equations for this $\Spc$ structure
with $\eta=0$. Gluing this solution to the solutions on $M$ given by the invariant $SW_-(M,\fs_M)$, we find 
$SW_-(X,\fs')=\pm 1$, where $\fs'$ is the $\Spc$ structure on $X$ obtained from $\fs_M$ and $\fs_N'$, compare Proposition~2
of~\cite{KMT}. We have $\fs'=\fs_{K^{-1}}\otimes E$, with $E=a_1e_1$. Therefore $E^2=-a_1^2\leq -1$, and
Corollary~\ref{minimal} shows that $X$ is not minimal. 
This completes the proof of Theorem~\ref{sum}.
\end{proof}

\bibliographystyle{amsplain}

\begin{thebibliography}{10}

\bibitem{BPV}
W.~Barth, C.~Peters and A.~Van~de~Ven, {\sl Compact Complex Surfaces},
Springer-Verlag, Berlin 1984.

\bibitem{B}
N.~Buchdahl, {\em On compact K\"ahler surfaces}, Ann.~Inst.~Fourier {\bf 49} (1999),
287--302.

\bibitem{D}
S.~K.~Donaldson, {\em The orientation of Yang--Mills moduli spaces and $4$-manifold topology}, 
J.~Differential Geometry {\bf 26} (1987), 397--428.

\bibitem{Go}
R.~E.~Gompf, {\em A new construction of symplectic manifolds}, 
Ann.~Math.~{\bf 142} (1995), 527--595.

\bibitem{GoS}
R.~E.~Gompf, {\em The topology of symplectic manifolds}, 
Turkish J.~Math.~{\bf 25} (2001), 43--59.

\bibitem{tani}
D.~Kotschick, {\em Four-manifold invariants of finitely presentable groups},
in {\sl Topology, Geometry and Field Theory}, ed.~K.~Fukaya et.~al.,
World Scientific 1994.

\bibitem{irred}
D.~Kotschick, {\em On irreducible four--manifolds}, Preprint 
alg-geom/9504012 (unpublished).

\bibitem{Bourbaki}
D.~Kotschick, {\em The Seiberg-Witten invariants of symplectic 
four--manifolds}, S\'eminaire Bourbaki, 48\`eme ann\'ee, 1995-96, 
no.~812, Ast\'erisque {\bf 241} (1997), 195--220.

\bibitem{KMT}
D.~Kotschick, J.~W.~Morgan and C.~H.~Taubes, {\em Four-manifolds 
without symplectic structures but with non-trivial Seiberg--Witten 
invariants}, Math. Research Letters {\bf 2} (1995), 119--124.

\bibitem{LL2} T.-J.~Li and A.-K.~Liu, {\em The equivalence between SW and Gr in 
the case where $b^+=1$}, Intern.~Math.~Res.~Notices {\bf 7} (1999), 335--345.

\bibitem{LL3} T.-J.~Li and A.-K.~Liu, {\em Uniqueness of symplectic canonical 
class, surface cone and symplectic cone of $4$-manifolds with $b^+=1$}, 
J.~Differential Geometry {\bf 58} (2001), 331--370.

\bibitem{Liu} A.-K.~Liu, {\em Some new applications of general wall crossing 
formula, Gompf's conjecture and its applications}, 
Math.~Research Letters {\bf 3} (1996), 569--585.

\bibitem{McD1}
D.~McDuff, {\em The structure of rational and ruled symplectic $4$-manifolds},
Journ.~Amer.~Math.~Soc.~{\bf 3} (1990), 679--712; Erratum: 
Journ.~Amer.~Math.~Soc.~{\bf 5} (1992), 987--988.

\bibitem{McD2}
D.~McDuff, {\em Blow ups and symplectic embeddings in dimension $4$},
Topology {\bf 30} (1991), 409--421. 

\bibitem{McD3}
D.~McDuff, {\em Immersed spheres in symplectic $4$-manifolds},
Ann.~Inst.~Fourier {\bf 42} (1992), 369--392. 

\bibitem{MP}
D.~McDuff and L.~Polterovich, {\em Symplectic packings and algebraic geometry},
Invent.~math.~{\bf 115} (1994), 405--429.

\bibitem{MS}
D.~McDuff and D.~Salamon, {\sl $J$-holomorphic Curves and Quantum 
Cohomology}, American Math.~Soc., Providence, R.~I.~1994.

\bibitem{Rjdg}
Y.~Ruan, {\em Symplectic topology on algebraic $3$-folds},  J.~Differential Geom.~{\bf 39} (1994), 
215--227.

\bibitem{R}
Y.~Ruan, {\em Symplectic topology and complex surfaces}, in {\sl Geometry and 
Analysis on Complex Manifolds}, Festschrift for Prof.~S.~Kobayashi's 60th Birthday,
ed.~T.~Mabuchi, J.~Noguchi and T.~Ochiai, World Scientific, Singapore 1994.

\bibitem{T1} 
C.~H.~Taubes, {\em The Seiberg--Witten invariants
and symplectic forms}, Math.~Research Letters {\bf 1} (1994),
809--822.

\bibitem{swg}
C.~H.~Taubes, {\em The Seiberg--Witten and the Gromov invariants},
Math.~Research Letters {\bf 2} (1995), 221--238.

\bibitem{main}
C.~H.~Taubes, {\em SW $\Rightarrow$ Gr, From the Seiberg--Witten
equations to pseudo--holomorphic curves}, Journal Amer.~Math.~Soc.~{\bf 
9} (1996), 845--918.

\bibitem{Taubes}
C.~H.~Taubes, {\sl Seiberg Witten and Gromov invariants for 
symplectic $4$-manifolds}, International Press 2000.

\end{thebibliography}

\bigskip

\end{document}